\documentclass[11pt]{article}
\usepackage{amsfonts}
\usepackage{color}
\newtheorem{theo}{Theorem}[section]
\newtheorem{lem}[theo]{Lemma}
\newtheorem{cor}[theo]{Corollary}

\newtheorem{defi}{Definition}[section]

\newcommand{\mysection}[1]{\section{#1} \setcounter{equation}{0}}
\newcommand{\proof}{{\sc Proof.} \quad}
\newcommand{\proofc}{{\sc Proof} \ }
\newcommand{\be}{\begin{equation} \label}
\newcommand{\ee}{\end{equation}}
\newcommand{\bea}{\begin{eqnarray}\label}
\newcommand{\eea}{\end{eqnarray}}
\newcommand{\bas}{\begin{eqnarray*}}
\newcommand{\eas}{\end{eqnarray*}}
\newcommand{\bit}{\begin{itemize}}
\newcommand{\eit}{\end{itemize}}
\newcommand{\qed}{\hfill$\Box$ \vskip.2cm}
\newcommand{\nn}{\nonumber}
\newcommand{\R}{\mathbb{R}}
\newcommand{\N}{\mathbb{N}}
\newcommand{\pO}{\partial\Omega}

\newcommand{\eps}{\varepsilon}

\newcommand{\supp}{{\rm supp} \, }

\newcommand{\wto}{\rightharpoonup}
\newcommand{\wsto}{\stackrel{\star}{\rightharpoonup}}
\newcommand{\hra}{\hookrightarrow}

\newcommand{\abs}{\\[5pt]}

\newcommand{\io}{\int_\Omega}

\newcommand{\ue}{u_\eps}
\newcommand{\ve}{v_\eps}
\newcommand{\we}{w_\eps}

\begin{document}
\title{Large-data global generalized solutions in a chemotaxis system with tensor-valued sensitivities}
\author{
Michael Winkler\footnote{michael.winkler@math.uni-paderborn.de}\\
{\small Institut f\"ur Mathematik, Universit\"at Paderborn,}\\
{\small 33098 Paderborn, Germany} }
\date{}
\maketitle
\begin{abstract}
\noindent 
  The chemotaxis system
  \bas
	\left\{ \begin{array}{l}
	u_t=\Delta u - \nabla \cdot \Big(uS(x,u,v)\cdot\nabla v\Big),
	\qquad x\in \Omega, \ t>0, \\[1mm]
	v_t=\Delta v - uf(v),
	\qquad x\in \Omega, \ t>0, \\[1mm]
 	\end{array} \right.
	\qquad \qquad (\star)
  \eas
  for the density $u=u(x,t)$ of a cell population and the concentration $v=v(x,t)$ of an attractive  
  chemical consumed by the former,
  is considered under no-flux boundary conditions in a bounded domain $\Omega\subset\R^n$, $n\ge 1$, 
  with smooth boundary, where 
  $f \in C^1([0,\infty);[0,\infty))$ and $S \in C^2(\bar\Omega\times [0,\infty)^2;\R^{n\times n})$ are given functions
  such that $f(0)=0$.\\
  In contrast to related Keller-Segel-type problems with scalar sensitivities, 
  in presence of such matrix-valued $S$ the system 
  ($\star$) in general apparently does not possess any useful gradient-like structure.
  Accordingly, its analysis needs to be based on new types of a priori bounds. \\
  Using a spatio-temporal $L^2$ estimate for $\nabla \ln (u+1)$ as a starting point,
  we derive a series of compactness properties of solutions to suitably regularized versions of ($\star$).
  Motivated by these, we develop a generalized solution concept 
  which requires solutions to satisfy very mild regularity hypotheses only, especially for the component $u$; 
  in particular, the chemotactic flux $uS(x,u,v)\cdot\nabla v$ needs not be integrable in this context.\\
  On the basis of the above compactness properties, it is finally shown that 
  within this framework, under a mild growth assumption on $S$ and for all sufficiently regular
  nonnegative initial data, the corresponding initial-boundary value problem for ($\star$) possesses at least one 
  global generalized solution.
  This extends known results which in the case of such general matrix-valued $S$ provide statements on global existence
  only in the two-dimensional setting and under the additional restriction that $\|v_0\|_{L^\infty(\Omega)}$ be small.\abs
\end{abstract}
\newpage
\mysection{Introduction}
{\bf Chemotaxis with tensor-valued sensitivities.}\quad
This work is concerned with solutions of the parabolic initial-boundary value problem
\be{0}
	\left\{ \begin{array}{l}
	u_t=\Delta u - \nabla \cdot \Big(uS(x,u,v)\cdot\nabla v\Big),
	\qquad x\in \Omega, \ t>0, \\[1mm]
	v_t=\Delta v - uf(v),
	\qquad x\in \Omega, \ t>0, \\[1mm]
	\nabla u \cdot \nu = u(S(x,u,v) \cdot \nabla v) \cdot \nu, \quad \nabla v \cdot \nu=0, 
	\qquad x\in \partial\Omega, \ t>0, \\[1mm]
	u(x,0)=u_0(x), \quad v(x,0)=v_0(x),
	\qquad x\in\Omega,
 	\end{array} \right.
\ee
in a bounded domain $\Omega \subset \R^n$ with smooth boundary, where $n\ge 1$ and $\nu$ denotes the outward
normal vector field on $\pO$, and where $f:[0,\infty)\to \R$ and the matrix-valued function 
$S:\Omega\times [0,\infty)^2 \to \R^{n\times n}$ are supposed
to be given parameter functions.\abs
Systems of this type arise in mathematical biology as models for the evolution of cell populations, in which
individuals, besides moving randomly, are able to partially adapt their motion to gradients of a chemical signal substance.
This mechanism, also known as {\em chemotaxis}, 
in prototypical situations is such that the preferred direction of motion 
is either toward increasing signal concentrations, or away from the latter (\cite{hillen_painter2009}).
A simple model for these processes of {\em chemoattractive} and {\em chemorepulsive} movement 
was proposed by Keller and Segel in 1970 (\cite{KS}), at its core containing an equation of the form
\be{ks1}
	u_t=\Delta u - \chi \nabla \cdot (u\nabla v)
\ee
for the evolution of the population density $u=u(x,t)$ in response to the gradient of the chemical concentration
$v=v(x,t)$, where the constant $\chi\in\R$ is positive in the attractive and negative in the repulsive case.
Such Keller-Segel-type systems, obtained upon 
complementing (\ref{ks1}) or variants thereof by appropriate equations for the chemical in the respective situations,
have widely been used as models in quite diverse particular biological contexts, including 
spontaneous aggregation phenomena in populations of {\em Dictyostelium discoideum} (\cite{KS}),
tumor cell invasion (\cite{chaplain_lolas}), and also self-organization during embryonic development 
(\cite{painter_maini_othmer2000}).\abs
%
%
%
%
%
In contrast to this, more recent modeling approaches (\cite{xue_othmer}, \cite{othmer_hillen2002}, \cite{diluzio_et_al})
suggest to allow for more general mechanisms of chemotactic
migration in certain situations, including directions not necessarily parallel to the gradient of the signal. 
Corresponding models then require the so-called chemotactic sensitivity, in (\ref{ks1}) represented by the constant scalar
$\chi$, to be a general matrix-valued function such as in (\ref{0}). 
For instance, a concise derivation of a macroscopic model for the behavior of swimming bacteria near the surface of their 
surrounding fluid, as presented in \cite{xue_othmer},
in the corresponding parabolic limit leads to a description of the cell population density by the first 
equation in (\ref{0}), with sensitivity tensors of the form
\be{S_rot}
	S(x,u,v)= 
	\chi \cdot \left( \begin{array}{cc}
	1 & 0 \\ 0 & 1
	\end{array} \right)
	+ \beta \cdot \left( \begin{array}{cc}
	0 & 1 \\ -1 & 0
	\end{array} \right),
\ee
in their nondiagonal parts inter alia reflecting that when cells swim e.g.~parallel to a surface, larger viscous forces 
are exerted on those parts of the cells which are closer to the surface.
Here in the simplest conceivable setting $\chi$ and $\beta$ are assumed to be positive constants, but they may as well vary
with $x$ such as e.g.~in cases when rotational flux components are neglected far from boundary regions,
and moreover possibly depend on the 
the variables $u$ and $v$ if further mechanisms are accounted for such as saturation effects 
at large cell or signal densities
(\cite{xue2013}, \cite{hillen_painter2009}).\abs
{\bf Boundedness vs.~blow-up.}\quad
Guided by this example, in this work we will concentrate on the case when besides such a general type of chemotactic 
motion, the coupling between the quantities $u$ and $v$ is governed by signal consumption through cells; 
that is, we shall assume that cells absorb the chemical in question upon contact, 
as reflected in the particular form of the second equation in (\ref{0}). 
A fundamental mathematical question is then whether and in which sense the resulting initial-boundary value problem (\ref{0})
can be solved globally in time. 
Since in view of the choice of boundary conditions, the system (\ref{0}) formally preserves the total mass of cells in the
sense that $\io u(\cdot,t)\equiv \io u_0$, addressing this question essentially amounts to either ruling out
or showing the occurrence of finite-time mass accumulation. 
Indeed, in the setting of the minimal version of the full original Keller-Segel system,
\be{ks}
	\left\{ \begin{array}{l}
	u_t=\Delta u - \nabla \cdot (u\nabla v), 
	\qquad x\in\Omega, \ t>0, \\[1mm]
	v_t=\Delta v - v + u,
	\qquad x\in\Omega, \ t>0, 
	\end{array} \right.
\ee
in which the signal is thus {\em produced} by the cells, 
such a singularity formation, mathematically represented as finite-time blow-up of the solution component $u$
with respect to the norm in $L^\infty(\Omega)$, may occur in certain situations:
For appropriate initial data, explosions of this type have been detected
when either $\Omega\subset \R^2$ is a disk and the total mass of cells is supercritical in the sense that $\io u_0>8\pi$
(\cite{herrero_velazquez}, \cite{mizoguchi_win}), or when $\Omega$ is a ball in $\R^n$, $n\ge 3$, and $\io u_0$ is an
arbitrary prescribed number (\cite{win_JMPA}).
Here the criticality of the spatially two-dimensional setting is underlined by a complementing result which asserts that
in this case the condition $\io u_0<4\pi$ is sufficient to ensure global existence of bounded solutions, thereby
ruling out any blow-up phenomenon; in the radial case, this condition can even be relaxed to the essentially optimal
inequality $\io u_0<8\pi$ (\cite{nagai_senba_yoshida}).
In the one-dimemsional version of (\ref{ks}), all solutions are global and bounded (\cite{osaki_yagi}), whereas
in the three- and higher-dimensional case alternative smallness assumptions on the initial data, involving 
norms of $u_0$ in $L^p(\Omega)$ for $p\ge \frac{n}{2}$, warrant global boundedness (\cite{win_JDE}, \cite{cao}).\abs
On the other hand, signal consumption as in (\ref{0}) is known to inhibit this tendency toward blow-up, at least to a
certain extent, when coupled to the mechanisms in (\ref{ks1}): For instance, the corresponding Neumann initial-boundary
value problem for the prototypical system
\be{ks_abs}
	\left\{ \begin{array}{l}
	u_t=\Delta u - \nabla \cdot (u\nabla v), 
	\qquad x\in\Omega, \ t>0, \\[1mm]
	v_t=\Delta v - uv,
	\qquad x\in\Omega, \ t>0, 
	\end{array} \right.
\ee
possesses global classical solutions in smoothly bounded convex domains $\Omega\subset\R^2$ for all reasonably regular
initial data, and moreover all these solutions approach the constant equilibrium given by $u\equiv \frac{1}{|\Omega|} \io u_0$
and $v\equiv 0$ in the large time limit (\cite{taowin_jde_absorb}). 
For the three-dimensional analogue, at least certain generalized global solutions can be constructed. 
These eventually become bounded and smooth and stabilize in the aforementioned manner, but it is unknown whether
they may develop singularities at an intermediate stage (\cite{taowin_jde_absorb}).
Related systems involving nonlinear cell diffusion, essentially modeled by a porous medium-type operator
$\Delta u^m$ or non-degenerate variants thereof, have recently been studied in \cite{cao_ishida} 
and \cite{wang_mu_zhou_ZAMP}, 
where it has been shown that global bounded solutions can be constructed
if the enhancement of diffusion at high densities is sufficiently large in the sense 
that $m>2-\frac{2}{n}$.\\
Some of these global existence and 
boundedness properties of (\ref{ks_abs}) can even be found in a more complex model for swimming
aerobic bacteria which, in addition to the mechanisms reflected in (\ref{ks_abs}), includes the interaction of 
cells and chemoattractant with the surrounding fluid (\cite{duan_lorz_markowich}, \cite{win_CPDE_fluid}, \cite{win_ARMA}).\abs
{\bf The mathematical challenge: Deriving boundedness despite loss of energy structure.}\quad
From a point of view of mathematical analysis, passing from (\ref{ks_abs}) to (\ref{0}) by 
allowing for more complex cross-diffusion mechanisms in (\ref{0}) appears to bring about a significant structural change:
For (\ref{ks_abs}), the integral
\bas
	\io \Big\{ u\ln u + 2 |\nabla \sqrt{v}|^2 \Big\}
\eas
plays the role of an energy functional in that it decreases along trajectories (\cite{taowin_jde_absorb}, cf.~also
\cite{duan_lorz_markowich}).
A corresponding gradient-like structure, along with all its consequences for the a priori knowledge on the regularity of 
solutions, apparently cannot be expected for general matrix-valued sensitivities $S$ in (\ref{0}).
It is thus not clear how far the blow-up preventing effect of signal absorption in (\ref{ks_abs}) extends
to the general system (\ref{0}).
As far as we know, the only available result in this direction asserts
global existence of bounded solutions to (\ref{0}) in bounded convex planar domains, 
even in the classical sense, under mild assumptions on $S$ and $f$
(essentially coinciding with (\ref{reg_f}), (\ref{reg_S}) and (\ref{S}) below), but only 
under the restrictive additional assumption that $\|v_0\|_{L^\infty(\Omega)}$ be small enough (\cite{lsxw}). 
Without such a smallness condition, the recent paper \cite{cao_ishida} proves 
global existence of bounded weak solutions to the related system obtained from (\ref{0})
upon replacing $\Delta u$ by the porous medium-type nonlinear diffusion term $\Delta u^m$ with arbitrary $m>1$.\abs
{\bf Main results.} \quad
The purpose of the present paper is to establish a result on global existence for (\ref{0}) under fairly general
assumptions on $f$ and $S$.		
More precisely, throughout our analysis we will assume that
\be{reg_f}
	f \in C^1([0,\infty)) \quad \mbox{is nonnegative with} \quad
	f(0)=0,
\ee
and that $S=(S_{ij})_{i,j\in \{1,...,n\}}$ is a chemotactic sensitivity tensor with
\be{reg_S}
	S_{ij} \in C^2(\bar\Omega \times [0,\infty) \times [0,\infty)) \qquad \mbox{for } i,j\in\{1,...,n\}.
\ee
Moreover, we suppose that with some nondecreasing function $S_0$ on $[0,\infty)$,
$S$ satisfies the growth hypothesis
\be{S}
	|S(x,u,v)| \le S_0(v)
	\qquad \mbox{for all } (x,u,v) \in \bar\Omega \times [0,\infty) \times [0,\infty).
\ee
Our main result then says that within this framework, for all suitably smooth initial data the problem
(\ref{0}) is globally solvable in an appropriate generalized sense. 
In particular, unlike in \cite{lsxw} we do neither need to impose any smallness asumption on the initial data here,
nor do we require any restriction on the spatial dimension.
\begin{theo}\label{theo20}
  Suppose that $n\ge 1$ and that $\Omega\subset\R^n$ is a bounded domain with smooth boundary, and 
  let $f$ and $S$ satisfy (\ref{reg_f}), (\ref{reg_S}) and (\ref{S}).
  Then for any choice of nonnegative functions $u_0\in C^0(\bar\Omega)$ and $v_0 \in W^{1,\infty}(\Omega)$,
  the problem (\ref{0}) posseses at least one global generalized solution $(u,v)$ in the sense of Definition
  \ref{defi14}. 
  This solution can be obtained as the limit a.e.~in $\Omega\times (0,\infty)$ 
  of a sequence $((\ue,\ve))_{\eps=\eps_j\searrow 0}$ 
  of smooth classical solutions to the regularized problems (\ref{0eps}) below.
\end{theo}
{\bf Key steps in our analysis.} \quad
In order to highlight the main ideas underlying our approach, and to outline the structure of this work,
let us note that unlike the case when $n=2$ and $\|v_0\|_{L^\infty(\Omega)}$ is assumed to be small enough,	
a priori estimates for the solution component $u$ in some reflexive Lebesgue space seem hard to obtain. 
As seen in \cite{lsxw} for convex planar domains,
such an additional smallness assumption indeed allows for the derivation of bounds for both $u$ and $|\nabla v|^2$
in $L^\infty((0,T);L^p(\Omega))$ with arbitrarily large $p>1$ through an essentially straightforward approach
using suitable differential inequalities for $\io u^p + \io |\nabla v|^{2p}$.
Instead, our analysis needs to be based on alternative a priori information on solutions $(\ue,\ve)$
to adequately regularized versions of (\ref{0}) (see (\ref{0eps}) below).
Here beyond the immediate boundedness properties associated with the conservation of mass functional $\io \ue(\cdot,t)$
and the nonincrease of $\|\ve(\cdot,t)\|_{L^\infty(\Omega)}$ (Lemma \ref{lem1}),
of fundamental importance to our approach will be the key estimate 
\be{key}
	\int_0^\infty \io \frac{|\nabla\ue|^2}{(\ue+1)^2} \le C
\ee
with some $C>0$ independent of the regularization parameter $\eps\in (0,1)$ (Lemma \ref{lem3}).
Due to the strong dampening at large values of $\ue$ of the weight function $\frac{1}{(\ue+1)^2}$ therein, however,
we do not expect (\ref{key}) to initiate an appropriate bootstrap process yielding substantial 
further regularity properties which would allow for passing to the limit $\eps\searrow 0$ suitably so as to
obtain a limit object solving (\ref{0}) in one of the standard weak formulations.
We shall accordingly introduce a generalized solution concept, to be specified in Definitions \ref{defi15},
\ref{defi14} and \ref{defi16}, which at its core refers to the transformed quantity $\ln (u+1)$ rather than to 
$u$ itself.\\
Indeed, viewing (\ref{key}) as an inequality for $\nabla \ln (\ue+1)$ and establishing an appropriate estimate for
$\partial_t \ln (\ue+1)$, we will thereby infer in Corollary \ref{cor7} that
\be{lnu}
	(\ln (\ue+1))_{\eps\in (0,1)}
	\quad \mbox{is relatively compact in $L^2_{loc}(\bar\Omega\times [0,\infty))$ with respect to the strong topology.}
\ee
Furthermore, (\ref{key}) will be essential in deriving in Lemma \ref{lem10} that 
\be{uf}
	(\ue f(\ve))_{\eps\in (0,1)}
	\quad \mbox{is relatively compact in $L^1_{loc}(\bar\Omega\times [0,\infty))$ 
	with respect to the weak topology.}
\ee
This will on the one hand allow for passing to the limit along suitable subsequences so as to obtain
a limit object $(u,v)$, for which $v$ solves the second equation in (\ref{0}) in the natural weak sense
(Section \ref{sect_passing}).
On the other hand, (\ref{uf}) will enable us to refine 
straightforward compactness properties of $(\ve)_{\eps\in (0,1)}$, as expressed in Section \ref{sect_comp_v}
and in Section \ref{sect_passing}, so as to obtain in Section \ref{sect_strong_nablav} that	
\be{nablav_strong}
	(\nabla \ve)_{\eps\in (0,1)}
	\quad \mbox{is relatively compact in $L^2_{loc}(\bar\Omega\times [0,\infty))$ with respect to the strong topology.}
\ee
In the natural weak version of the first equation in (\ref{0eps}) associated with $\ln (\ue+1)$ (see (\ref{17.1})),
these compactness properties (\ref{lnu}) and (\ref{nablav_strong}) will form a main ingredient in 
taking $\eps\searrow 0$ termwise,	
with one exception being an integral containing $\frac{|\nabla u|^2}{(u+1)^2}$, 
for which it seems that only a one-sided control can be obtained 
by using lower semicontinuity of norms with respect to weak convergence.\\
Therefore, in our solution concept we shall require $\ln (u+1)$ to satisfy the respective
integral {\em inequality} only, thus
generalizing a {\em supersolution} property of $u$ with regard to the first equation in (\ref{0}) (Definition \ref{defi14}).
As seen in Lemma \ref{lem5}, the role of a complementing {\em subsolution}-like property can be played by the 
simple nonincrease of mass (cf.~(\ref{16.1})),
and within this framework the above limit $(u,v)$ indeed is a generalized solution of (\ref{0})
(see Section \ref{sect_solution_u}).
\mysection{A generalized solution concept}\label{sect_defi}
To begin with, let us first specify our solution concept. 
As far as the second component $v$ is concerned, a generalization of the respective sub-problem of (\ref{0})
is rather straightforward, because there the only nonlinear part $uf(v)$ is of lowest differentiability order.
\begin{defi}\label{defi15}
  Let $u\in L^1_{loc}(\bar\Omega \times [0,\infty))$, let $f$ satisfy (\ref{reg_f}), and assume that $v_0\in W^{1,2}(\Omega)$.
  Then a nonnegative function
  \bas
	v\in L^\infty_{loc}(\bar\Omega \times [0,\infty)) \cap L^2_{loc}([0,\infty);W^{1,2}(\Omega))
  \eas
   is said to be a {\em global weak solution} of
  \be{15.1}
	\left\{ \begin{array}{l}
	v_t=\Delta v - uf(v), \qquad x\in\Omega, \ t>0, \\[1mm]
	\frac{\partial v}{\partial\nu}=0, \qquad x\in \pO, \ t>0,\\[1mm]
	v(x,0)=v_0(x), \qquad x\in \Omega,
	\end{array} \right.
  \ee
  if for all $\varphi\in L^\infty(\Omega \times (0,\infty)) \cap L^2((0,\infty);W^{1,2}(\Omega))$ having compact support
  in $\bar\Omega \times [0,\infty)$ with $\varphi_t\in L^2(\Omega\times (0,\infty))$, the identity
  \be{15.2}
	\int_0^\infty \io v\varphi_t + \io v_0\varphi(\cdot,0)
	= \int_0^\infty \io \nabla v \cdot \nabla \varphi + \int_0^\infty \io uf(v)\varphi
  \ee
  holds.
\end{defi}
The most important part of our solution concept refers to the cross-diffusive equation in (\ref{0}).
\begin{defi}\label{defi14}
  Assume that $S$ complies with (\ref{reg_S}) and (\ref{S}), and that $u_0\in L^\infty(\Omega)$ is nonnegative.
  Moreover, let $\phi \in C^2([0,\infty))$ be such that $\phi'>0$ on $(0,\infty)$, and suppose that 
  $v\in L^\infty_{loc}(\bar\Omega \times [0,\infty)) \cap L^2_{loc}([0,\infty);W^{1,2}(\Omega))$ is nonnegative.
  Then a nonnegative function $u:\Omega \times (0,\infty) \to \R$ will be called a {\em global very weak
  $\phi$-supersolution} of the problem
  \be{14.1}
	\left\{ \begin{array}{l}
	u_t=\Delta u - \nabla \cdot \Big( uS(x,u,v)\cdot \nabla v \Big),
	\qquad x\in \Omega, t>0, \\[1mm]
	\Big(\nabla u - u(S(x,u,v)\cdot \nabla v) \Big)\cdot \nu=0, 
	\qquad x\in\pO, \ t>0, \\[1mm]
	u(x,0)=u_0(x), \qquad x\in\Omega,
	\end{array} \right.
  \ee
  if
  \bea{14.2}
	& & \mbox{$\phi(u)$ and $\phi''(u)|\nabla u|^2$ belong to $L^1_{loc}(\bar\Omega \times [0,\infty))$,} \nn\\[1mm]
	& & \mbox{$u\phi''(u)\nabla u$ and $u\phi'(u)$ belong to $L^2_{loc}(\bar\Omega \times [0,\infty))$,}
  \eea
  and if for each nonnegative $\varphi \in C_0^\infty(\bar\Omega \times [0,\infty))$
  with $\frac{\partial\varphi}{\partial\nu}=0$ on $\pO \times (0,\infty)$, the inequality
  \bea{14.3}
	- \int_0^\infty \io \phi(u) \varphi_t - \io \phi(u_0) \varphi(\cdot,0)
	&\ge& \int_0^\infty \io \phi(u) \Delta\varphi
	- \int_0^\infty \io \phi''(u) |\nabla u|^2 \varphi \nn\\
	& &
	+ \int_0^\infty \io u\phi''(u) \nabla u \cdot \Big(S(x,u,v) \cdot \nabla v \Big) \varphi \nn\\
	& &
	+ \int_0^\infty \io u\phi'(u) \Big(S(x,u,v)\cdot\nabla v \Big) \cdot\nabla \varphi
  \eea
  is satisfied.
\end{defi}
{\bf Remark.} \quad
  i) \ It can easily be checked using (\ref{S}) that the required regularity properties of $v$ along with (\ref{14.2})
  ensure that all integrals in (\ref{14.3}) are well-defined.\\
  ii) \ In our final existence argument given in Lemma \ref{lem17}, we shall eventually choose $\phi(s):=\ln (s+1)$
  for $s\ge 0$.\abs
It is evident that in order to become meaningful, 
the above supersolution property has to be complemented by an additional condition which
rules out that the component $u$ has its time derivative exceeding the one dictated by the first equation in (\ref{0}).
We shall see that in a generalized sense, for this it is already sufficient to require that only the total mass 
$\io u(\cdot,t)$ be bounded {\em from above} by $\io u_0$:
\begin{defi}\label{defi16}
  A couple $(u,v)$ of nonnegative functions $u$ and $v$ defined in $\Omega \times (0,\infty)$ and satisfying
  \bas
	u \in L^\infty((0,\infty);L^1(\Omega))
  \eas
  with
  \be{16.1}
	\io u(\cdot,t) \le \io u_0
	\qquad \mbox{for a.e.~$t>0$}
  \ee
  as well as
  \bas
	v \in L^\infty_{loc}(\bar\Omega \times [0,\infty))
	\cap L^2_{loc}([0,\infty);W^{1,2}(\Omega))
  \eas
  will be named a {\em global generalized solution} of (\ref{0}) if $v$ is a global weak solution of (\ref{15.1})
  in the sense of Definition \ref{defi15},
  and if for some $\phi \in C^2([0,\infty))$ with $\phi'>0$ on $[0,\infty)$, $u$ is a global very weak $\phi$-supersolution
  of (\ref{14.1}) in the sense of Definition \ref{defi14}.
\end{defi}
Indeed, this concept is fully compatible with that of classical solutions in the following sense:
\begin{lem}\label{lem5}
  Suppose that $u$ and $v$ are nonnegative functions from $C^0(\bar\Omega \times [0,\infty)) \cap
  C^{2,1}(\bar\Omega \times (0,\infty))$. 
  Then if $(u,v)$ is a global generalized solution of (\ref{0}), it follows that $(u,v)$ also is a classical solution
  of (\ref{0}) in $\Omega \times (0,\infty)$.
\end{lem}
\proof
  Since it is clear upon a standard reasoning that $v$ is a classical solution of (\ref{15.1}), we only need to prove
  that $u$ is a classical solution of (\ref{14.1}).
  To this end, we first fix a sequence $(\zeta_j)_{j\in\N} \subset C_0^\infty([0,\infty))$ such that
  $0 \le \zeta_j \le 1=\zeta(0)$, $\zeta_j' \le 0$ and $\supp \zeta_j \subset [0,\frac{1}{j}]$ for $j\in\N$,
  and given any nonnegative $\psi \in C_0^\infty(\Omega)$ 
  we choose $\varphi(x,t):=\zeta_j(t)\psi(x)$, $(x,t)\in \bar\Omega\times [0,\infty)$,
  in (\ref{14.3}). Then thanks to (\ref{14.2}), the dominated convergence theorem and the fact that $\zeta_j'$
  approaches the Dirac measure $-\delta(t)$, in the limit $j\to\infty$ we obtain
  \bas
	\io \phi(u(\cdot,0))\psi - \io \phi(u_0) \psi \ge 0
  \eas
  for any such $\psi$. This implies that $\phi(u(\cdot,0)) \ge \phi(u_0)$ in $\Omega$ and hence $u(\cdot,0) \ge u_0$
  in $\Omega$, because $\phi'>0$ on $[0,\infty)$. Therefore, (\ref{16.1}) and the continuity of $u$ at $t=0$ warrant that
  actually $u(\cdot,0)=u_0$ in $\Omega$.\\
  Secondly, choosing
  arbitrary nonnegative $\varphi \in C_0^\infty(\Omega \times (0,\infty))$ in (\ref{14.3}), by a similar density argument
  we see that
  \bas
	\frac{\partial}{\partial t} \phi(u)
	\ge \Delta \phi(u) - \phi''(u) |\nabla u|^2
	- \phi'(u) \nabla \cdot \Big(uS(x,u,v)\cdot \nabla v\Big)
	\qquad \mbox{in } \Omega \times (0,\infty)
  \eas
  holds in the classical sense. This is equivalent to
  \bas
	\phi'(u) u_t \ge \phi'(u) \Delta u - \phi'(u) \nabla \cdot \Big(uS(x,u,v)\cdot \nabla v\Big)
	\qquad \mbox{in } \Omega \times (0,\infty),
  \eas
  and using that $\phi'>0$ on $[0,\infty)$ we conclude that $u$ is a classical supersolution of the first equation in
  (\ref{0}), that is,
  \be{5.1}
	u_t \ge \Delta u - \nabla \cdot \Big(uS(x,u,v)\cdot \nabla v\Big)
	\qquad \mbox{in } \Omega \times (0,\infty).
  \ee
  Finally, choosing arbitrary nonnegative $\varphi \in C_0^\infty(\bar\Omega \times (0,\infty))$ supported near 
  $\pO \times [0,\infty)$ and such that $\frac{\partial\varphi}{\partial\nu}=0$ on $\pO$, in a standard manner we 
  moreover obtain from (\ref{14.3}) that
  \be{5.2}
	\frac{\partial u}{\partial\nu} \ge u \Big(S(x,u,v)\cdot\nabla v\Big) \cdot \nu
	\qquad \mbox{on } \pO \times (0,\infty).
  \ee
  Now if $u$ was not a classical solution of (\ref{14.1}) then, by (\ref{5.1}), (\ref{5.2}) and a continuity argument,
  for some open subset $G_1 \subset \Omega$ and some open interval $J_1\subset (0,\infty)$ we would have
  \be{5.3}
	u_t>\Delta u - \nabla \cdot \Big(uS(x,u,v)\cdot\nabla v \Big)
	\qquad \mbox{in } G_1 \times J_1,
  \ee
  or there would exist a relatively open set $G_2\subset \pO$ and an open interval $J_2\subset (0,\infty)$ fulfilling
  \be{5.4}
	\frac{\partial u}{\partial\nu} > u\Big(S(x,u,v)\cdot\nabla v\Big)
	\qquad \mbox{in } G_2 \times J_2.
  \ee
  In the former case, (\ref{5.3}) together with (\ref{5.1}) and (\ref{5.2}) would imply that for all $t\in J_1$,
  \bas
	\io u(\cdot,t) - \io u_0 = \int_0^t \io u_t
	&>& \int_0^t \io \bigg(\Delta u - \nabla \cdot \Big(uS(x,u,v)\cdot\nabla v\Big) \bigg) \\
	&=& \int_0^t \int_{\pO} \bigg(\frac{\partial u}{\partial\nu} - \Big(uS(x,u,v)\cdot\nabla v\Big) \cdot \nu \bigg) 
	\\[2mm]
	&\ge& 0,
  \eas
  meaning that $\io u(\cdot,t) > \io u_0$ for all $t\in J_1$ and thereby contradicting the second assumption (\ref{16.1})
  on $u$. Along with a similar argument in the case when (\ref{5.4}) holds, this completes the proof.
\qed
\mysection{Global solutions of regularized problems}
In order to introduce an appropriate regularization of (\ref{0}), let us fix families $(\rho_\eps)_{\eps\in (0,1)}$
and $(\chi_\eps)_{\eps\in (0,1)}$ of functions 
\bas
	\rho_\eps \in C_0^\infty(\Omega)
	\quad \mbox{such that} \quad
	0 \le \rho_\eps \le 1 \mbox{ in $\Omega$ \quad and \quad}
	\rho_\eps \nearrow 1 \mbox{ in $\Omega$ as $\eps\searrow 0$}
\eas
and
\bas
	\chi_\eps \in C_0^\infty([0,\infty))
	\quad \mbox{such that} \quad
	0 \le \chi_\eps \le 1 \mbox{ in $[0,\infty)$ \quad and \quad}
	\chi_\eps \nearrow 1 \mbox{ in $[0,\infty)$ as $\eps\searrow 0$,}
\eas
define
\bas
	S_\eps(x,u,v):=\rho_\eps(x) \cdot \chi_\eps(u) \cdot S(x,u,v),
	\qquad x\in\bar\Omega, \ u\ge 0, \ v \ge 0,
\eas
and consider the problems
\be{0eps}
	\left\{ \begin{array}{l}
	u_{\eps t} = \Delta \ue - \nabla \cdot \Big(\ue S_{\eps}(x,\ue,\ve) \cdot \nabla \ve \Big),
	\qquad x\in \Omega, \ t>0, \\[1mm]
	v_{\eps t} = \Delta \ve - \ue f(\ve),
	\qquad x\in \Omega, \ t>0, \\[1mm]
	\frac{\partial \ue}{\partial\nu} = \frac{\partial \ve}{\partial\nu}=0
	\qquad x\in \partial\Omega, \ t>0, \\[1mm]
	\ue(x,0)=u_0(x), \quad \ve(x,0)=v_0(x), 
	\qquad x\in\Omega,
 	\end{array} \right.
\ee
for $\eps\in (0,1)$. 
These are indeed globally solvable in the classical sense:
\begin{lem}\label{lem0}
  For all $\eps\in (0,1)$, there exists a pair
  $(\ue,\ve) \in (C^0(\bar\Omega \times [0,\infty)) \cap C^{2,1}(\bar\Omega \times (0,\infty))$ of nonnegative functions
  which solve (\ref{0eps}) classically in $\Omega \times (0,\infty)$.
\end{lem}
\proof
  Local existence of a smooth solution can be seen by a well-established contraction mapping
  argument in the space $C^0(\bar\Omega \times [0,T]) \times L^\infty((0,T);W^{1,q}(\Omega))$ 
  for arbitrary fixed $q>\max\{2,n\}$
  and suitably small $T>0$ (see \cite{win_cpde}, for instance). 
  Since $S_\eps(x,u,v)\equiv 0$ for all sufficiently large $u$, standard estimation techniques yield extensibility of 
  this local solution for all times (cf.~e.g.~\cite{horstmann_win}).
\qed
The following basic properties of solutions to (\ref{0eps}) are immediate.
\begin{lem}\label{lem1}
  The solution of (\ref{0eps}) satisfies
  \be{mass}
	\io \ue(\cdot,t)=\io u_0 
	\qquad \mbox{for all } t>0
  \ee
  as well as
  \be{vinfty}
	\|\ve(\cdot,t)\|_{L^\infty(\Omega)} \le \|v_0\|_{L^\infty(\Omega)}
	\qquad \mbox{for all } t>0.
  \ee
  In particular, with $S_0$ as defined in (\ref{S}) we have the pointwise estimate
  \be{1.3}
	|S_\eps(x,\ue,\ve)| \le S_1:=S_0(\|v_0\|_{L^\infty(\Omega)})
	\qquad \mbox{in } \Omega \times (0,\infty).
  \ee
\end{lem}
\proof
  The identity (\ref{mass}) directly results upon integration of the first equation in (\ref{0eps}) with respect to
  $x\in\Omega$. The estimate (\ref{vinfty}) is a straightforward consequence of the maximum principle applied
  to the second equation in (\ref{0eps}), because we already know that $\ue \ge 0$, and because $f$ was assumed to be
  nonnegative throughout.
\qed
Two more testing procedures easily yield further information:
\begin{lem}\label{lem2}
  The solution of (\ref{0eps}) has the properties
  \be{nablav}
	\int_0^\infty \io |\nabla \ve|^2 \le \frac{1}{2} \io v_0^2
  \ee
  and
  \be{2.2}
	\int_0^\infty \io \ue f(\ve) \le \io v_0.
  \ee
\end{lem}
\proof
  Multiplying the second equation in (\ref{0eps}) by $\ve$ and integrating by parts over $\Omega$, we obtain
  \bas
	\frac{1}{2} \frac{d}{dt} \io \ve^2 + \io |\nabla \ve|^2 = - \io \ue \ve f(\ve)
	\qquad \mbox{for all } t>0.
  \eas
  Since here by nonnegativity of $f$, $\ue$ and $\ve$ the right-hand side is nonpositive, integrating in time 
  yields (\ref{nablav}).\\
  Likewise, testing the second equation in (\ref{0eps}) against a nontrivial constant shows that
  \bas
	\frac{d}{dt} \io \ve = - \io \ue f(\ve)
	\qquad \mbox{for all } t>0,
  \eas
  from which (\ref{2.2}) results upon a time integration.
\qed
\mysection{Estimates for $\ln(\ue+1)$}
We proceed to derive further estimates for $\ue$.
The first of these provides an integral bound for the gradient of $\ln (\ue+1)$.
\begin{lem}\label{lem3}
  For each $\eps\in (0,1)$, the solution of (\ref{0eps}) satisfies
  \be{3.1}
	\int_0^\infty \io \frac{|\nabla \ue|^2}{(\ue+1)^2} \le K_1:=2\io u_0 + \frac{S_1^2}{2} \cdot \io v_0^2,
  \ee
  with the number $S_1$ being as defined in (\ref{1.3}).
\end{lem}
\proof
  We multiply the first equation in (\ref{0eps}) by $\frac{1}{\ue+1}$ and integrate by parts over $\Omega$, which results in
  the identity
  \bea{3.2}
	\frac{d}{dt} \io \ln (\ue+1) 
	&=& - \io \nabla \ue \cdot \nabla \frac{1}{\ue+1}
	+ \io \nabla \frac{1}{\ue+1} \cdot \Big( \ue S_\eps(x,\ue,\ve) \cdot \nabla \ve \Big) \nn\\
	&=& \io \frac{|\nabla \ue|^2}{(\ue+1)^2}
	- \io \frac{\ue}{(\ue+1)^2} \nabla \ue \cdot \Big( S_\eps(x,\ue,\ve) \cdot \nabla \ve \Big)
	\quad \mbox{for all } t>0.
  \eea
  By Young's inequality and (\ref{1.3}),
  \bas
	\bigg| - \io \frac{\ue}{(\ue+1)^2} \nabla \ue \cdot \Big( S_\eps(x,\ue,\ve) \cdot \nabla \ve \Big) \bigg|
	&\le& \frac{1}{2} \io \frac{|\nabla \ue|^2}{(\ue+1)^2}\\
	& & + \frac{1}{2} \io \frac{\ue^2}{(\ue+1)^2} \cdot |S_\eps(x,\ue,\ve)|^2 \cdot |\nabla \ve|^2 \\
	&\le& \frac{1}{2} \io \frac{|\nabla \ue|^2}{(\ue+1)^2}
	+ \frac{S_1^2}{2} \cdot \io |\nabla \ve|^2.
  \eas
  Therefore, an integration of (\ref{3.2}) with respect to the time variable yields
  \bas
	\io \ln (\ue(\cdot,t)+1) - \io \ln (u_0+1)
	&\ge& \frac{1}{2} \int_0^t \io \frac{|\nabla \ue|^2}{(\ue+1)^2}
	- \frac{S_1^2}{2} \int_0^t \io |\nabla \ve|^2
	\qquad \mbox{for all } t>0
  \eas
  and thus, since $0\le \ln (\xi+1) \le \xi$ for all $\xi\ge 0$,
  \bas
	\frac{1}{2} \int_0^t \io \frac{|\nabla \ue|^2}{(\ue+1)^2}
	&\le& \io \ue(\cdot,t)
	+\frac{S_1^2}{2} \int_0^t \io |\nabla \ve|^2 \\
	&=& \io u_0 
	+\frac{S_1^2}{2} \int_0^t \io |\nabla \ve|^2
	\qquad \mbox{for all } t>0
  \eas
  thanks to (\ref{mass}). As $\int_0^t \io |\nabla \ve|^2 \le \frac{1}{2} \io v_0^2$ by Lemma \ref{lem2}, this
  establishes (\ref{3.1}).
\qed
In order to prepare pointwise convergence a.e.~in $\Omega\times (0,\infty)$ for $u_\eps$ along a suitable sequence
of numbers $\eps=\eps_j\searrow 0$, we next aim at deriving a strong compactness property of $(\ln (\ue+1))_{\eps \in (0,1)}$.
This is prepared by the following.
\begin{lem}\label{lem6}
  Let $m\in\N$ be such that $m>\frac{n}{2}$.
  Then there exists $K_2>0$ with the property that for each $\eps \in (0,1)$, the solution of (\ref{0eps}) satisfies
  \be{6.1}
	\int_0^T \Big\|\partial_t \ln (\ue (\cdot,t)+1)\Big\|_{(W_0^{m,2}(\Omega))^\star} dt \le K_2 \cdot (1+T)
	\qquad \mbox{for all } T>0.
  \ee
\end{lem}
\proof
  For fixed $t>0$ and arbitrary $\psi \in W_0^{m,2}(\Omega)$, using the first equation in (\ref{0eps}) and integrating 
  by parts we obtain
  \bea{6.2}
	\io \partial_t \ln (\ue(x,t)+1) \cdot \psi(x)dx
	&=& \io \frac{u_{\eps t}}{\ue+1} \cdot \psi \nn\\
	&=& \io \frac{1}{\ue+1} \Delta \ue \cdot \psi
	- \io \frac{1}{\ue+1} \nabla \cdot \Big( \ue S_\eps(x,\ue,\ve)\cdot\nabla \ve \Big) \psi \nn\\
	&=& - \io \frac{1}{\ue+1} \nabla \ue \cdot \nabla \psi
	+ \io \frac{1}{(\ue+1)^2} |\nabla \ue|^2 \psi \nn\\
	& & + \io \frac{\ue}{\ue+1} \Big(S_\eps(x,\ue,\ve)\cdot \nabla \ve\Big) \cdot\nabla \psi \nn\\
	& & - \io \frac{\ue}{(\ue+1)^2} \nabla \ue \cdot \Big(S_\eps(x,\ue,\ve)\cdot\nabla \ve\Big) \psi.
  \eea
  Here, by the Cauchy-Schwarz inequality we have
  \bas
	\bigg| - \io \frac{1}{\ue+1} \nabla \ue \cdot \nabla \psi \bigg| 
	\le \bigg( \io \frac{|\nabla \ue|^2}{(\ue+1)^2} \bigg)^\frac{1}{2} \cdot \|\nabla\psi\|_{L^2(\Omega)},
  \eas
  and by the same token we see that
  \bas
	\bigg| \io \frac{\ue}{\ue+1} \Big(S_\eps(x,\ue,\ve)\cdot\nabla \ve \Big)\cdot\nabla \psi \bigg|
	\le S_1 \cdot \bigg(\io |\nabla \ve|^2 \bigg)^\frac{1}{2} \cdot \|\nabla \psi\|_{L^2(\Omega)}
  \eas
  and
  \bas
	\bigg| 	- \io \frac{\ue}{(\ue+1)^2} \nabla \ue \cdot \Big(S_\eps(x,\ue,\ve)\cdot\nabla \ve\Big) \psi \bigg|
	\le S_1 \cdot \bigg( \io \frac{|\nabla \ue|^2}{(\ue+1)^2} \bigg)^\frac{1}{2} 
	\cdot \bigg(\io |\nabla \ve|^2 \bigg)^\frac{1}{2} \cdot \|\psi\|_{L^\infty(\Omega)}.
  \eas
  Since clearly 
  \bas
	\bigg| \io \frac{1}{(\ue+1)^2} |\nabla \ue|^2 \psi \bigg|
	\le \bigg(\io \frac{|\nabla \ue|^2}{(\ue+1)^2} \bigg) \cdot \|\psi\|_{L^\infty(\Omega)},
  \eas
  (\ref{6.2}) therefore yields
  \bas
	\bigg| \io \partial_t \ln (\ue(x,t)+1) \cdot \psi(x)dx \bigg|
	&\le& \Bigg\{ \bigg(\io \frac{|\nabla \ue|^2}{(\ue+1)^2} \bigg)^\frac{1}{2}
	+ \io \frac{|\nabla\ue|^2}{(\ue+1)^2} \\
	& & \hspace*{5mm}
	+ S_1 \cdot \bigg(\io |\nabla \ve|^2 \bigg)^\frac{1}{2}
	+ S_1 \cdot \bigg(\io \frac{|\nabla \ue|^2}{(\ue+1)^2} \bigg)^\frac{1}{2} 
		\cdot \bigg(\io |\nabla\ve|^2 \bigg)^\frac{1}{2} \Bigg\} \times \\
	& & \times \bigg( \|\nabla \psi\|_{L^2(\Omega)} + \|\psi\|_{L^\infty(\Omega)} \bigg)
	\qquad \mbox{for all } \psi \in W_0^{m,2}(\Omega).
  \eas
  As our condition $m>\frac{n}{2}$ ensures that 
  the space $W_0^{m,2}(\Omega)$ is continuously embedded into $L^\infty(\Omega)$,
  by Young's inequality this implies that with some $c_1>0$,
  \bas
	\bigg| \io \partial_t \ln (\ue(x,t)+1) \cdot \psi(x)dx \bigg|
	\le c_1 \cdot \bigg\{ 1 + \io \frac{|\nabla\ue|^2}{(\ue+1)^2} + \io |\nabla \ve|^2 \bigg\} \cdot
	\|\psi\|_{W_0^{m,2}(\Omega)}
  \eas
  for all $\psi\in W_0^{n,2}(\Omega)$, meaning that
  \bas
	\Big\|\partial_t \ln(\ue(\cdot,t)+1) \Big\|_{(W_0^{m,2}(\Omega))^\star} 
	\le c_1 \cdot \bigg\{ 1 + \io \frac{|\nabla\ue|^2}{(\ue+1)^2} + \io |\nabla \ve|^2 \bigg\} 
	\qquad \mbox{for all } t>0.
  \eas
  Since according to Lemma \ref{lem3} and (\ref{nablav}) we have $\int_0^\infty \io \frac{|\nabla\ue|^2}{(\ue+1)^2} \le K_1$
  and $\int_0^\infty \io |\nabla \ve|^2 \le \frac{1}{2} \io v_0^2$, an integration over $(0,T)$ easily yields
  (\ref{6.1}) with an evident choice of $K_2$.
\qed
Now a straightforward application of (a variant of) the Aubin-Lions lemma can be used to establish 
the following compactness properties of $(\ln (\ue+1))_{\eps\in (0,1)}$.
\begin{cor}\label{cor7}
  Let $T>0$. Then $(\ln (\ue+1))_{\eps\in (0,1)}$ is relatively compact in $L^2((0,T);W^{1,2}(\Omega))$ 
  with respect to the weak topology, and relatively compact in $L^2(\Omega \times (0,T))$ with respect to the
  strong topology.
\end{cor}
\proof
  As $(\ln(\ue+1))_{\eps\in (0,1)}$ is bounded in $L^2((0,T);W^{1,2}(\Omega))$ according to Lemma \ref{lem3} and
  (\ref{mass}), the first statement is immediate. 
  Using that moreover $(\partial_t \ln (\ue+1))_{\eps\in (0,1)}$ is bounded in $L^1((0,T);(W_0^{n,2}(\Omega))^\star)$
  by Lemma \ref{lem6}, since $(W_0^{n,2}(\Omega))^\star$ is a Hilbert space we may invoke a version of the 
  Aubin-Lions lemma (\cite[Theorem 2.3]{temam}) to obtain the claimed strong precompactness property.
\qed
\mysection{Compactness properties of $(\ve)_{\eps\in (0,1)}$}\label{sect_comp_v}
By a simplified variant of the argument of the previous section, we can readily derive the following.
\begin{lem}\label{lem9}
  Let $T>0$. Then $(\ve)_{\eps\in (0,1)}$ is relatively compact in $L^2(\Omega\times (0,T))$ with respect to the strong
  topology.
\end{lem}
\proof
  We let $m\in\N$ be such that $m>\frac{n}{2}$, and take an arbitrary $\psi \in W_0^{m,2}(\Omega)$. 
  Then from the second equation in (\ref{0eps})
  and the Cauchy-Schwarz inequality, for each fixed $t\in (0,T)$ we obtain
  \bas
	\bigg| \io v_{\eps t}(x,t)\psi(x)dx \bigg|
	&=& \bigg| - \io \nabla \ve \cdot \nabla \psi - \io \ue f(\ve) \psi \bigg| \\
	&\le& \bigg(\io |\nabla \ve|^2 \bigg)^\frac{1}{2} \cdot \|\nabla\psi\|_{L^2(\Omega)}
	+ \bigg( \io \ue f(\ve) \bigg) \cdot \|\psi\|_{L^\infty(\Omega)}.
  \eas
  Again since $W_0^{m,2}(\Omega) \hra L^\infty(\Omega)$, we thus find that
  \bas
	\int_0^T \|v_{\eps t}(\cdot,t)\|_{(W_0^{m,2}(\Omega))^\star} dt
	\le c_1 \int_0^T \bigg\{ 1 + \io |\nabla \ve|^2 + \io \ue f(\ve) \bigg\} dt
  \eas
  with some $c_1>0$, and hence in light of Lemma \ref{lem2} we conclude that
  \bas
	\int_0^T \|v_{\eps t}(\cdot,t)\|_{(W_0^{m,2}(\Omega))^\star} dt
	\le c_1 T + \frac{c_1}{2} \io v_0^2 + c_1 \io v_0.
  \eas
  Therefore, the Aubin-Lions lemma in \cite[Theorem 2.3]{temam} along with the boundedness of $(\ve)_{\eps\in (0,1)}$ in 
  $L^2((0,T);W^{1,2}(\Omega))$, as asserted by (\ref{vinfty}) and (\ref{nablav}), yields the claim.
\qed
\mysection{Precompactness of $(\ue f(\ve))_{\eps\in (0,1)}$}
In passing to the limit in the taxis term in (\ref{0eps}), we will also need strong precompactness
of $(\nabla\ve)_{\eps\in (0,1)}$ in $L^2_{loc}(\bar\Omega \times [0,\infty))$, rather than the corresponding 
weak compactness property implied by (\ref{nablav}). 
This will finally be achieved in Lemma \ref{lem13} below, but prepared by a series of steps, 
the first of which can be interpreted as providing some superlinear integrability
property of the inhomogeneity $h_\eps:=\ue f(\ve)$ in the semilinear heat equation $v_{\eps t}=\Delta v_\eps - h_\eps$.
\begin{lem}\label{lem4}
  For each $\eps \in (0,1)$ we have the inequality
  \be{4.1}
	\int_0^\infty \io \ue \ln (\ue+1) f(\ve)
	\le K_3,
  \ee
  where
  \bas
	K_3:=\io v_0 \ln (u_0+1) + (\|v_0\|_{L^\infty(\Omega)} + 2)\cdot K_1
	+ \Big(\frac{1}{2} + \frac{S_1}{2} + \frac{1}{8} \|v_0\|_{L^\infty(\Omega)}^2 S_1^2 \Big) \cdot \io v_0^2
  \eas
  with $S_1$ and $K_1$ taken from (\ref{1.3}) and (\ref{3.1}), respectively.
\end{lem}
\proof
  Using the first and second equation in (\ref{0eps}), we track the time evolution of $\io \ve \ln (\ue+1)$ by computing
  \bea{4.2}
	\frac{d}{dt} \io \ve \ln (\ue+1)
	&=& \io v_{\eps t} \ln (\ue+1) 
	+\io \frac{\ve}{\ue+1} u_{\eps t} \nn\\
	&=& \io \Delta \ve \cdot \ln (\ue+1)
	- \io \ue \ln (\ue+1) f(\ve) \nn\\
	& & + \io \frac{\ve}{\ue+1} \Delta \ue
	- \io \frac{\ve}{\ue+1} \nabla \cdot \Big( \ue S_\eps(x,\ue,\ve) \cdot \nabla \ve \Big)
  \eea 
   for all $t>0$. Integrating by parts, we find that
  \bas
	\io \Delta \ve \cdot \ln (\ue+1)
	= - \io \frac{1}{\ue+1} \nabla \ue \cdot \nabla \ve
  \eas
  and
  \bas
	\io \frac{\ve}{\ue+1} \Delta \ue
	&=& - \io \nabla \Big(\frac{\ve}{\ue+1}\Big) \cdot \nabla \ue \\
	&=& - \io \frac{1}{\ue+1} \nabla \ue \cdot \nabla \ve
	+ \io \frac{\ve}{(\ue+1)^2} |\nabla \ue|^2
  \eas
  as well as
  \bas
	- \io \frac{\ve}{\ue+1} \nabla \cdot \Big( \ue S_\eps(x,\ue,\ve) \cdot \nabla \ve \Big)
	&=& \io \ue \nabla \Big(\frac{\ve}{\ue+1}\Big) \cdot \Big(S_\eps(x,\ue,\ve) \cdot \nabla \ve \Big) \\
	&=& \io \frac{\ue}{\ue+1} \nabla \ve \cdot \Big(S_\eps(x,\ue,\ve)\cdot \nabla \ve \Big) \\
	& & - \io \frac{\ue\ve}{(\ue+1)^2} \nabla \ue \cdot \Big(S_\eps(x,\ue,\ve)\cdot \nabla \ve \Big)
  \eas
  for $t>0$.
  Upon a time integration, (\ref{4.2}) therefore becomes
  \bea{4.3}
	\int_0^t \io \ue \ln (\ue+1) f(\ve)
	+ \io \ve(\cdot,t) \ln (\ue(\cdot,t)+1)
	&=& \io v_0 \ln (u_0+1) \nn\\
	& & \hspace*{-15mm}
	- 2\int_0^t \io \frac{1}{\ue+1} \nabla \ue \cdot \nabla \ve
	+ \int_0^t \io \frac{\ve}{(\ue+1)^2} |\nabla \ue|^2 \nn\\
	& & \hspace*{-15mm}
	+ \int_0^t \io \frac{\ue}{\ue+1} \nabla \ve \cdot \Big(S_\eps(x,\ue,\ve) \cdot \nabla \ve\Big)\nn\\
	& & \hspace*{-15mm}
	- \int_0^t \io \frac{\ue\ve}{(\ue+1)^2} \nabla \ue \cdot \Big(S_\eps(x,\ue,\ve) \cdot \nabla \ve \Big)
  \eea
  for all $t>0$.
  Here we use Young's inequality, Lemma \ref{lem3} and (\ref{nablav}) in estimating
  \bas
	- 2\int_0^t \io \frac{1}{\ue+1} \nabla \ue \cdot \nabla \ve
	&\le& \int_0^t \io \frac{|\nabla \ue|^2}{(\ue+1)^2} 
	+ \int_0^t \io |\nabla \ve|^2 \\
	&\le& K_1 + \frac{1}{2} \io v_0^2
	\qquad \mbox{for all } t>0,
  \eas
  whereas Lemma \ref{lem3} combined with (\ref{vinfty}) shows that
  \bas
	\int_0^t \io \frac{\ve}{(\ue+1)^2} |\nabla \ue|^2 
	&\le& \|v_0\|_{L^\infty(\Omega)} \cdot \int_0^t \io \frac{|\nabla\ue|^2}{(\ue+1)^2}
	\le \|v_0\|_{L^\infty(\Omega)} \cdot K_1
	\qquad \mbox{for all } t>0.
  \eas
  Moreover, by means of (\ref{1.3}) and (\ref{nablav}) we find that
  \bas
	\int_0^t \io \frac{\ue}{\ue+1} \nabla \ve \cdot \Big(S_\eps(x,\ue,\ve \cdot \nabla \ve\Big)
	&\le& S_1 \cdot \int_0^t \io |\nabla \ve|^2 \\
	&\le& S_1 \cdot \frac{1}{2} \io v_0^2
	\qquad \mbox{for all } t>0,
  \eas
  and similarly (\ref{3.1}), (\ref{vinfty}) and (\ref{nablav}) in view of Young's inequality yield
  \bas
	- \int_0^t \io \frac{\ue\ve}{(\ue+1)^2} \nabla \ue \cdot \Big(S_\eps(x,\ue,\ve) \cdot \nabla \ve \Big)
	&\le& \int_0^t \io \frac{|\nabla \ue|^2}{(\ue+1)^2} \\
	& & + \frac{1}{4} \|v_0\|_{L^\infty(\Omega)}^2 \cdot 
		S_1^2 \cdot \int_0^t \io \frac{\ue^2}{(\ue+1)^2} |\nabla \ve|^2 \\
	&\le& K_1 + \frac{1}{4} \|v_0\|_{L^\infty(\Omega)}^2 \cdot S_1^2 \cdot \frac{1}{2} \io v_0^2
  \eas
  for all $t>0$.
  Since $\io \ve \ln (\ue+1)$ is nonnegative, (\ref{4.3}) therefore implies (\ref{4.1}).
\qed
Along with the Pettis theorem, the above lemma yields the following.
\begin{lem}\label{lem10}
  For each $T>0$, the family $(\ue f(\ve))_{\eps \in (0,1)}$ is relatively compact in $L^1(\Omega \times (0,T))$ with respect
  to the weak topology.
\end{lem}
\proof
  Let $\we:=\ue f(\ve)$, $\eps\in (0,1)$. Then since
  $f(\ve) \le c_1:=\|f\|_{L^\infty((0,\|v_0\|_{L^\infty(\Omega)}))}$ according to (\ref{vinfty}), using 
  Lemma \ref{lem2} and Lemma \ref{lem4}
  and writing $c_2:=\max \{1,c_1\}$ and $m:=\io u_0$ we find that
  \bas
	\int_0^T \io \we \ln (\we+1)
	&\le& \int_0^T \io \ue f(\ve) \cdot \ln \big( c_1 \ue + 1 \big) \\
	&\le& \int_0^T \io \ue f(\ve) \cdot \ln \Big(c_2 (\ue+1)\Big) \\
	&=& \ln c_2 \cdot \int_0^T \io \ue f(\ve)
	+ \int_0^T \io \ue \ln (\ue+1) f(\ve) \\[2mm]
	&\le& \ln c_2 \cdot \io v_0 + K_3.
  \eas
  In view of Pettis' theorem, this equi-integrability property already guarantees that $(\we)_{\eps\in (0,1)}$
  is relatively compact with respect to the weak topology in $L^1(\Omega \times (0,T))$.
\qed
\mysection{Passing to the limit. Solution properties of $v$}\label{sect_passing}
We can now perform a first subsequence extraction procedure, resulting in a limit object $(u,v)$ the second component
of which can already be shown to be a weak solution of its respective equation in (\ref{0}).
\begin{lem}\label{lem11}
  There exists a sequence $(\eps_j)_{j\in\N}$ of numbers $\eps_j \in (0,1)$ such that $\eps_j \searrow 0$ as $j\to\infty$ and
  \begin{eqnarray}
	& & \ue \to u \qquad \mbox{a.e.~in } \Omega \times (0,\infty), \label{11.1} \\
	& & \ln (\ue+1) \wto \ln (u+1) \qquad \mbox{in } L^2_{loc}([0,\infty);W^{1,2}(\Omega)), \label{11.2} \\
	& & \ve \to v \qquad \mbox{a.e.~in } \Omega \times (0,\infty), \label{11.3} \\
	& & \ve \to v \qquad \mbox{in } L^2_{loc}(\bar\Omega \times [0,\infty)), \label{11.4} \\
	& & \ve \wsto v \qquad \mbox{in } L^\infty(\Omega \times (0,\infty)), \label{11.5} \\
	& & \nabla \ve \wto \nabla v \qquad \mbox{in } L^2(\Omega \times (0,\infty)) \qquad \qquad \mbox{and} \label{11.55} \\
	& & \ue f(\ve) \to uf(v) \qquad \mbox{in } L^1_{loc}(\bar\Omega \times [0,\infty)) \label{11.6}
  \eea
  as $\eps=\eps_j\searrow 0$ with certain nonnegative functions $u$ and $v$ defined in $\Omega \times (0,\infty)$.
  Moreover, $v$ is a weak solution of (\ref{15.1}) in the sense of Definition \ref{defi15}.
\end{lem}
\proof
  According to (\ref{vinfty}), (\ref{nablav}) and Lemma \ref{lem9}, (\ref{11.3})-(\ref{11.55}) can be achieved through
  a straightforward extraction process. 
  Similarly, Corollary \ref{cor7} and Lemma \ref{lem3} imply that (\ref{11.1}) and (\ref{11.2}) hold along a further
  subsequence. 
  In particular, by continuity of $f$ this entails that
  \be{11.8}
	\ue f(\ve) \to uf(v) \qquad \mbox{a.e.~in } \Omega\times (0,\infty)
  \ee
  as $\eps=\eps_j\searrow 0$, which combined with Lemma \ref{lem10} and Egorov's theorem ensures that upon another
  extraction we may assume that
  \bas
	\ue f(\ve) \wto uf(v)
	\qquad \mbox{in } L^1_{loc}(\bar\Omega \times [0,\infty))
  \eas
  as $\eps=\eps_j\searrow 0$.
  In light of Lemma \ref{lem8} below, again using (\ref{11.8}) we conclude that even (\ref{11.6}) holds.\\
  Now the verification of the claimed solution property of $v$ is quite standard: Given $\varphi$ with the properties
  listed in Definition \ref{defi15}, testing the second equation in (\ref{0eps}) against $\varphi$ yields
  \be{11.9}
	\int_0^\infty \io \ve \varphi_t + \io v_0 \varphi(\cdot,0)
	= \int_0^\infty \io \nabla \ve \cdot \nabla \varphi
	+ \int_0^\infty \io \ue f(\ve) \cdot \varphi
  \ee
  for all $\eps\in (0,1)$. Since $\varphi$ has compact support in $\bar\Omega \times [0,\infty)$, the properties
  $\varphi_t\in L^2(\Omega\times (0,\infty)), \nabla \varphi \in L^2(\Omega\times (0,\infty))$ and 
  $\varphi \in L^\infty(\Omega \times (0,\infty))$ in conjunction with (\ref{11.4}), (\ref{11.55}) and (\ref{11.6}),
  respectively, imply that the identity (\ref{15.2}) results from (\ref{11.9}) upon taking
  $\eps=\eps_j\searrow 0$ in each integral separately.
\qed
\mysection{Strong precompactness of $(\nabla \ve)_{\eps \in (0,1)}$}\label{sect_strong_nablav}
Let us next fully concentrate on the problem of asserting strong precompactness of $(\nabla\ve)_{\eps\in (0,1)}$.
Having $\nabla v$ as a candidate for the desired limit at hand now, and knowing that by the weak convergence statement
in (\ref{11.55}) we have
$\int_0^T \io |\nabla v|^2 \le \liminf_{\eps=\eps_j\searrow 0} \int_0^T \io |\nabla \ve|^2$ for $T>0$,
in order to show that actually $\nabla \ve \to \nabla v$ in $L^2(\Omega \times (0,T))$ it is sufficient to make
sure that $\int_0^T \io |\nabla v|^2$ satisfies a corresponding estimate from below. 
This will be a consequence of the following lemma which is concerned with the standard entropy identity
\bas
	\frac{1}{2} \frac{d}{dt} \io v^2 + \io |\nabla v|^2 = - \io uvf(v), \qquad t>0,
\eas
that clearly holds for smooth solutions of (\ref{15.1}), 
but which seems not to be extensible in a straightforward way 
to arbitrary weak solutions $v$ of (\ref{15.1}) with the function $u$ on the right-hand side only belonging
to the non-reflexive space $L^\infty((0,T);L^1(\Omega))$.
After all, upon a suitable choice of test functions in (\ref{15.2}) it is possible to derive
a corresponding inequality which will be sufficient for our purpose.
\begin{lem}\label{lem12}
  There exists a null set $N\subset (0,\infty)$ such that
  the limit functions $u$ and $v$ gained in Lemma \ref{lem11} satisfy the inequality
  \be{12.1}
	\frac{1}{2} \io v^2(\cdot,T) - \frac{1}{2} \io v_0^2 + \int_0^T \io |\nabla v|^2 \ge - \int_0^T \io uvf(v)
	\qquad \mbox{for all } T\in (0,\infty) \setminus N.
  \ee
\end{lem}
\proof
  Since $v\in L^\infty(\Omega \times (0,\infty))$, $z(t):=\io v^2(x,t)dx, \ t>0$, defines a function 
  $z\in L^1_{loc}([0,\infty))$. Therefore there exists a null set $N\subset (0,\infty)$ such that each $T\in (0,\infty)
  \setminus N$ is a Lebesgue point of $z$; in particular, 
  \be{12.2}
	\frac{1}{\delta} \int_T^{T+\delta} \io v^2(x,t)dxdt
	\to \io v^2(x,T)dx
	\quad \mbox{for all } T\in (0,\infty) \setminus N
	\qquad \mbox{as } \delta\searrow 0.
  \ee
  To see that (\ref{12.1}) holds with this choice of $N$, given any $T\in (0,\infty) \setminus N$
  and $\delta\in (0,1)$ we let
  \bas
	\zeta_\delta(t) := \left\{ \begin{array}{ll}
	1, \qquad & t\in [0,T], \\[1mm]
	1-\frac{t-T}{\delta}, \qquad & t\in (T,T+\delta), \\[1mm]
	0, & t\ge T,
	\end{array} \right.
  \eas
  and define
  \bas
	\tilde v_k(x,t):=\left\{ \begin{array}{ll}
	v(x,t), \qquad & (x,t) \in \Omega \times (0,\infty), \\[1mm]
	v_{0k}(x), & (x,t) \in \Omega \times (-1,0],
	\end{array} \right.
  \eas
  for $k\in\N$, where $(v_{0k})_{k\in\N} \subset C^1(\bar\Omega)$ is such that $v_{0k} \to v_0$ in $L^2(\Omega)$.
  Then for $\delta\in (0,1), k\in\N$ and $h\in (0,1)$ we introduce
  \bas
	\varphi(x,t):=\varphi_{\delta,k,h}(x,t):=\zeta_\delta(t) \cdot (A_h \tilde v_k)(x,t),
	\qquad (x,t)\in \Omega \times (0,\infty),
  \eas
  where the temporal average $A_h \tilde v_k$ is defined as
  \bas
	(A_h \tilde v_k)(x,t):=\frac{1}{h} \int_{t-h}^t \tilde v_k(x,s)ds, 
	\qquad (x,t) \in \Omega \times (0,\infty).
  \eas
  Since $v\in L^\infty(\Omega\times (0,\infty)) \cap L^2((0,\infty);W^{1,2}(\Omega))$ by Lemma \ref{lem11},
  it can easily be checked that also $\varphi$ belongs to 
  $L^\infty(\Omega\times (0,\infty)) \cap L^2((0,\infty);W^{1,2}(\Omega))$, and that in addition $\varphi$ is supported in
  $\bar\Omega \times [0,T+1]$ with
  \bas
	\varphi_t(x,t)=\zeta_\delta'(t) \cdot (A_h \tilde v_k)(x,t) 
	+ \zeta_\delta(t) \cdot \frac{1}{h} \Big(\tilde v_k(x,t)-\tilde v_k(x,t-h)\Big),
	(x,t) \in \Omega \times (0,\infty),
  \eas
  implying that $\varphi_t \in L^2(\Omega \times (0,\infty))$.
  We may therefore insert $\varphi$ into (\ref{15.2}) to obtain
  \bea{12.3}
	I_1(\delta,k,h)+I_2(\delta,k,h) &:=&
	\int_0^\infty \io \zeta_\delta(t) \nabla v(x,t) \cdot \nabla (A_h \tilde v_k)(x,t)dxdt \nn\\
	& & + \int_0^\infty \io \zeta_\delta(t) u(x,t) f(v(x,t)) \cdot (A_h \tilde v_k)(x,t) dxdt \nn\\
	&=& \int_0^\infty \io \zeta_\delta'(t) v(x,t) \cdot (A_h \tilde v_k)(x,t)dxdt \nn\\
	& & + \int_0^\infty \io \zeta_\delta(t) v(x,t) \cdot 
		\frac{1}{h} \Big(\tilde v_k(x,t)-\tilde v_k(x,t-h)\Big)dxdt \nn\\
	& & + \io v_0(x) v_{0k} (x) dx \nn\\[2mm]
	&=:& I_3(\delta,k,h)+I_4(\delta,k,h)+I_5(\delta,k,h),
  \eea
  where we have used that
  \bas
	\varphi(x,0)=\zeta_\delta(0) \cdot \frac{1}{h} \int_{-h}^0 \tilde v_k(x,s) ds= v_{0k}(x), 
	\qquad x\in\Omega,
  \eas
  by definition of $\zeta_\delta$ and $\tilde v_k$.
  Now since $v_{0k} \in C^1(\bar\Omega)$, it follows that $\nabla \tilde v_k \in L^2(\Omega \times (-1,T+1))$, so that
  Lemma \ref{lem_a2} a) below applies to yield
  \bas
	\nabla (A_h\tilde v_k) = A_h (\nabla \tilde v_k) \wto \nabla \tilde v_k=\nabla v
	\quad \mbox{in } L^2(\Omega \times (0,T+1))
	\qquad \mbox{as } h\searrow 0,
  \eas
  so that
  \be{12.4}
	I_1(\delta,k,h) \to \int_0^\infty \io \zeta_\delta(t) \cdot |\nabla v|^2(x,t) dxdt
	\qquad \mbox{as } h\searrow 0.
  \ee
  Similarly, the inclusion $\tilde v_k\in L^\infty(\Omega \times (-1,T+1))$ along with Lemma \ref{lem_a2} b) ensures that
  \be{12.44}
	A_h \tilde v_k \wsto \tilde v_k=v
	\quad \mbox{in } L^\infty(\Omega \times (0,T+1))
	\qquad \mbox{as } h\searrow 0,
  \ee
  and that hence
  \be{12.5}
	I_2(\delta,k,h) \to \int_0^\infty \io \zeta_\delta(t) u(x,t) v(x,t) f(v(x,t)) dxdt
	\qquad \mbox{as } h\searrow 0,
  \ee
  because $uf(v) \in L^1((\Omega\times (0,T+1))$.
  By (\ref{12.44}) we clearly also see that
  \be{12.6}
	I_3(\delta,k,h) \to \int_0^\infty \io \zeta_\delta'(t) v^2(x,t)dxdt
	\qquad \mbox{as } h\searrow 0.
  \ee
  In order to analyze the corresponding limit behavior of $I_4(\delta,k,h)$, we split this integral according to
  \bas
	I_4(\delta,k,h)
	= \frac{1}{h} \int_0^\infty \io \zeta_\delta(t) \tilde v_k^2(x,t)dxdt
	- \frac{1}{h} \int_0^\infty \io \zeta_\delta(t) \tilde v_k(v,t) \tilde v_k(x,t-h) dxdt
  \eas
  and estimate the second term on the right by means of Young's inequality to see that
  \bas
	\frac{1}{h} \int_0^\infty \io \zeta_\delta(t) \tilde v_k(v,t) \tilde v_k(x,t-h) dxdt
	&\le& \frac{1}{2h} \int_0^\infty \io \zeta_\delta(t) \tilde v_k^2(x,t)dxdt \\
	& & + \frac{1}{2h} \int_0^\infty \io \zeta_\delta(t) \tilde v_k^2(x,t-h)dxdt.
  \eas
  Thus, upon substituting $s=t-h$, we find that
  \bas
	I_4(\delta,k,h)
	&\ge& \frac{1}{2h} \int_0^\infty \io \zeta_\delta(t) v^2(x,t)dxdt
	- \frac{1}{2h} \int_0^\infty \io \zeta_\delta(t) \tilde v_k^2(x,t-h)dxdt \\
	&=&  \frac{1}{2h} \int_0^\infty \io \zeta_\delta(t) v^2(x,t)dxdt
	- \frac{1}{2h} \int_0^\infty \io \zeta_\delta(s+h) v^2(x,s) dxds \\
	& & - \frac{1}{2h} \int_0^h \io \zeta_\delta(t) v_{0k}^2 (x) dxdt \\
	&=& - \frac{1}{2} \int_0^\infty \io \frac{\zeta_\delta(t+h)-\zeta_\delta(t)}{h} \cdot v^2(x,t)dxdt
	- \frac{1}{2h} \int_0^h \io \zeta_\delta(t) v_{0k}^2 (x) dxdt,
  \eas
  again because $\tilde v_k(\cdot,t)=v_{0k}$ for $t\in (-1,0)$. 
  Here since $\zeta_\delta$ is continuous with $\zeta_\delta(0)=0$, we have
  \bas
	- \frac{1}{2h} \int_0^h \io \zeta_\delta(t) v_{0k}^2 (x) dxdt
	\to -\frac{1}{2} \io v_{0k}^2(x)dx
	\qquad \mbox{as } h\searrow 0,
  \eas
  whereas by the dominated convergence theorem we conclude that
  \bas
	- \frac{1}{2} \int_0^\infty \io \frac{\zeta_\delta(t+h)-\zeta_\delta(t)}{h} \cdot v^2(x,t)dxdt
	\to - \frac{1}{2} \int_0^\infty \io \zeta_\delta'(t) v^2(x,t)dxdt
	\qquad \mbox{as } h\searrow 0,
  \eas
  whence altogether we infer that
  \bas
	\liminf_{h\searrow 0} I_4(\delta,k,h)
	\ge -\frac{1}{2} \int_0^\infty \io \zeta_\delta'(t) v^2(x,t) dxdt - \frac{1}{2} \io v_{0k}^2(x)dx.
  \eas
  Therefore, taking $h\searrow 0$ and recalling (\ref{12.4}), (\ref{12.5}) and (\ref{12.6}), from (\ref{12.3}) 
  we obtain the inequality
  \bas
	\int_0^\infty \io \zeta_\delta(t) |\nabla v(x,t)|^2 dxdt
	&+& \int_0^\infty \io \zeta_\delta(t) u(x,t) v(x,t) f(v(x,t)) dxdt \\
	&\ge& \frac{1}{2} \int_0^\infty \io \zeta_\delta'(t) v^2(x,t)dxdt \\
	& & - \frac{1}{2} \io v_{0k}^2(x)dx
	+ \io v_0(x) v_{0k}(x)dx
  \eas
  for all $k\in\N$, in the limit $k\to\infty$ implying that
  \bea{12.7}
	\int_0^\infty \io \zeta_\delta(t) |\nabla v(x,t)|^2 dxdt
	&+& \int_0^\infty \io \zeta_\delta(t) u(x,t)v(x,t)f(v(x,t)) dxdt \nn\\
	&\ge& \frac{1}{2} \int_0^\infty \io \zeta_\delta'(t) v^2(x,t) dxdt
	+ \frac{1}{2} \io v_0^2(x)dx.
  \eea
  Now by definition of $\zeta_\delta$, the first term on the right satisfies
  \bas
	\frac{1}{2} \int_0^\infty \io \zeta_\delta'(t) v^2(x,t) dxdt
	&=& -\frac{1}{2\delta} \int_T^{T+\delta} \io v^2(x,t)dxdt \\
	&\to& - \frac{1}{2} \io v^2(x,T)dx
	\qquad \mbox{as } \delta\searrow 0
  \eas
  according to the Lebesgue point property of $T$.
  Applying the monotone convergence theorem to both integrals on the left of (\ref{12.7}), we thereupon readily 
  arrive at (\ref{12.1}).
\qed
We can now establish the desired strong convergence result. Besides on the above inequality (\ref{12.1}), 
its derivation essentially relies on the strong convergence statement in (\ref{11.6}).
\begin{lem}\label{lem13}
  Let $(\eps_j)_{j\in\N}$ be as provided by Lemma \ref{lem11}. 
  Then there exists a subsequence, again denoted by $(\eps_j)_{j\in\N}$, such that for each $T>0$ we have
  \be{13.1}
	\nabla \ve \to \nabla v
	\quad \mbox{in } L^2(\Omega \times (0,T))
	\qquad \mbox{as } \eps=\eps_j\searrow 0.
  \ee
\end{lem}
\proof
  Since we know from (\ref{11.4}) that $\ve \to v$ in $L^2_{loc}(\bar\Omega \times [0,\infty))$, upon passing to
  a subsequence if necessary we may assume that as $\eps=\eps_j\searrow 0$ we have
  \bas
	\io \ve^2(\cdot,T) \to \io v^2(\cdot,T)
	\qquad \mbox{for all } T \in (0,\infty) \setminus N_1
  \eas
  with some null set $N_1 \subset (0,\infty)$. 
  Taking $N \subset (0,\infty)$ as in Lemma \ref{lem12}, we then evidently only need to verify (\ref{13.1}) 
  for all $T\in (0,\infty) \setminus (N\cup N_1)$. 
  Given any such $T$, we apply (\ref{11.6}) to see that as $\eps=\eps_j\searrow 0$,
  \bas
	\ue f(\ve) \to uf(v)
	\qquad \mbox{in } L^1(\Omega \times (0,T)),
  \eas
  which thanks to the fact that
  \bas
	\ve \wsto v
	\qquad \mbox{in } L^\infty(\Omega \times (0,T))
  \eas
  by (\ref{11.5}) implies that
  \bas
	\int_0^T \io \ue \ve f(\ve) \to \int_0^T \io uvf(v).
  \eas
  Therefore Lemma \ref{lem12} says that due to our choice of $T$,
  \bas
	\int_0^T \io |\nabla v|^2
	&\ge& - \frac{1}{2} \io v^2(\cdot,T)
	+ \frac{1}{2} \io v_0^2 - \int_0^T \io uvf(v) \\
	&=& \lim_{\eps=\eps_j\searrow 0} \bigg\{ -\frac{1}{2} \io \ve^2(\cdot,T) + \frac{1}{2} \io v_0^2
	- \int_0^T \io \ue \ve f(\ve) \bigg\}.
  \eas
  Since testing the second equation in (\ref{0eps}) by $\ve$ yields
  \bas
	-\frac{1}{2} \io \ve^2(\cdot,T) + \frac{1}{2} \io v_0^2 -\int_0^T \io \ue\ve f(\ve)
	= \int_0^T \io |\nabla \ve|^2,
  \eas
  this entails that
  \be{13.3}
	\int_0^T \io |\nabla v|^2 \ge \liminf_{\eps=\eps_j\searrow 0} \int_0^T \io |\nabla \ve|^2 .
  \ee
  On the other hand, by lower semicontinuity of the norm in $L^2(\Omega \times (0,T))$ with respect to weak convergence,
  \be{13.4}
	\int_0^T \io |\nabla v|^2 \le \liminf_{\eps=\eps_j\searrow 0} \int_0^T \io |\nabla \ve|^2,
  \ee
  whence (\ref{13.1}) results from (\ref{13.3}) and (\ref{13.4}) together with (\ref{11.55}) upon a well-known argument.
\qed
\mysection{Solution properties of $u$. Proof of Theorem \ref{theo20}}\label{sect_solution_u}
We are now in the position to show that also the limit $u$ from Lemma \ref{lem11} solves its associated
subsystem in (\ref{0}) in the sense specified in Definition \ref{defi16}.
We first establish the mass inequality (\ref{16.1}).
\begin{lem}\label{lem18}
  The function $u$ gained in Lemma \ref{lem11} satisfies $u\in L^\infty((0,\infty);L^1(\Omega))$ and
  \bas
	\io u(\cdot,t) \le \io u_0
	\qquad \mbox{for a.e.~$t>0$}.
  \eas
\end{lem}
\proof
  Since according to (\ref{mass}) we have $\io \ue(\cdot,t)=\io u_0$ for all $t>0$ and each $\eps\in (0,1)$,
  both statements are consequences from (\ref{11.1}) and Fatou's lemma.
\qed
The derivation of a corresponding $\phi$-supersolution property of $u$ in the spirit of Definition \ref{defi14}
is more delicate and crucially involves Lemma \ref{lem13}.
\begin{lem}\label{lem17}
  Let $u$ and $v$ be as constructed in Lemma \ref{lem11}. Then $u$ is a global very weak $\phi$-supersolution of
  (\ref{14.1}) with
  \bas
	\phi(s):=\ln (s+1), \qquad s\ge 0,
  \eas
  in the sense of Definition \ref{defi14}.
\end{lem}
\proof
  Using (\ref{11.2}), we first see that $\phi(u)$ and $u\phi'(u)=\frac{u}{u+1}$ belong to $L^1_{loc}(\bar\Omega \times
  [0,\infty))$ and to $L^2_{loc}(\bar\Omega \times [0,\infty))$, respectively. 
  Moreover, (\ref{11.2}) guarantees that
  \bas
	\phi''(u) |\nabla u|^2 = -\frac{|\nabla u|^2}{(u+1)^2} = -|\nabla \ln (u+1)|^2
	\in L^1_{loc}(\bar\Omega\times [0,\infty)),
  \eas
  and that since
  \bas
	|u\phi''(u)\nabla u| = \frac{u|\nabla u|}{(u+1)^2} \le |\nabla \ln (u+1)|,
  \eas
  we also have $u\phi''(u)\nabla u \in L^2_{loc}(\bar\Omega \times [0,\infty))$.\\
  Now given any nonnegative $\varphi \in C_0^\infty(\bar\Omega \times [0,\infty))$ with $\frac{\partial\varphi}{\partial\nu}
  =0$ on $\pO \times (0,\infty)$, multiplying the first equation in (\ref{0eps}) by 
  $\phi'(\ue) \cdot \varphi=\frac{1}{\ue+1}\cdot \varphi$ and integrating by parts we derive the identity
  \bea{17.1}
	\int_0^\infty \io \frac{1}{(\ue+1)^2} |\nabla \ue|^2 \varphi
	&=& - \int_0^\infty \io \ln (\ue+1)\varphi_t
	-\io \ln (u_0+1)\varphi(\cdot,0) \nn\\
	& & - \int_0^\infty \io \ln (\ue+1)\Delta\varphi \nn\\
	& & + \int_0^\infty \io \frac{\ue}{(\ue+1)^2} \nabla \ue \cdot \Big(S_\eps(x,\ue,\ve)\cdot\nabla \ve \Big) 
		\cdot \varphi \nn\\
	& & - \int_0^\infty \io \frac{\ue}{\ue+1} \Big(S_\eps(x,\ue,\ve)\cdot\nabla \ve\Big) \cdot \nabla \varphi
  \eea
  for all $\eps \in (0,1)$. Here, thanks to (\ref{11.2}) we have
  \be{17.2}
	- \int_0^\infty \io \ln (\ue+1)\varphi_t 
	\to - \int_0^\infty \io \ln (u+1)\varphi_t
  \ee
  and
  \be{17.3}
	- \int_0^\infty \io \ln (\ue+1) \Delta \varphi 
	\to - \int_0^\infty \io \ln (u+1)\Delta \varphi
  \ee
  as $\eps=\eps_j\searrow 0$. 
  Furthermore, by definition of $S_\eps$ the statements (\ref{11.1}) and (\ref{11.3}) imply that
  as $\eps=\eps_j\searrow 0$ we have
  \bas
	\frac{\ue}{\ue+1} S_\eps(x,\ue,\ve) \to \frac{u}{u+1} S(x,u,v)
	\qquad \mbox{a.e.~in } \Omega \times (0,\infty).
  \eas
  In light of the uniform majorization
  \bas
	\bigg| \frac{\ue}{\ue+1} S_\eps(x,\ue,\ve) \bigg|
	\le S_1
	\quad \mbox{in } \Omega \times (0,\infty) 
	\qquad \mbox{for all } \eps\in (0,1),
  \eas
  as warranted by Lemma \ref{lem1}, and the fact that due to Lemma \ref{lem13} we have
  \bas
	\nabla \ve \to \nabla v
	\qquad \mbox{in } L^2_{loc}(\bar\Omega\times [0,\infty)),
  \eas
  according to Lemma \ref{lem19} this implies the strong convergence property
  \be{17.4}
	\frac{\ue}{\ue+1} \Big(S_\eps(\cdot,\ue,\ve)\cdot\nabla \ve \Big)
	\to \frac{u}{u+1} \Big(S(\cdot,u,v)\cdot\nabla v\Big)
	\qquad \mbox{in } L^2_{loc}(\bar\Omega\times [0,\infty))
  \ee
  as $\eps=\eps_j\searrow 0$.
  Together with (\ref{11.2}) this entails that
  \bea{17.5}
	\int_0^\infty \io \frac{\ue}{(\ue+1)^2} \nabla \ue \cdot \Big( S_\eps(x,\ue,\ve) \cdot \nabla \ve \Big) \cdot \varphi
	&=& \int_0^\infty \io \nabla \ln (\ue+1) \cdot \Big( \frac{\ue}{\ue+1} S_\eps(x,\ue,\ve)\cdot \nabla \ve \Big)
		\cdot \varphi \nn\\
	&\to& \int_0^\infty \io \nabla \ln (u+1) \cdot \Big( \frac{u}{u+1} S(x,u,v)\cdot\nabla v\Big) \cdot \varphi \nn\\
	&=& \int_0^\infty \io \frac{u}{(u+1)^2} \nabla u \cdot \Big(S(x,u,v)\cdot\nabla v\Big) \cdot \varphi
  \eea
  as $\eps=\eps_j\searrow 0$. 
  Moreover, (\ref{17.4}) clearly also guarantees that
  \be{17.6}
	- \int_0^\infty \io \frac{\ue}{\ue+1} \Big(S_\eps(x,\ue,\ve)\cdot\nabla \ve \Big) \cdot\nabla \varphi
	\to -\int_0^\infty \io \frac{u}{u+1} \Big(S(x,u,v)\cdot\nabla v\Big) \cdot\nabla \varphi
  \ee
  as $\eps=\eps_j\searrow 0$. \\
  Collecting (\ref{17.2}), (\ref{17.3}), (\ref{17.5}) and (\ref{17.6}), by a lower semicontinuity argument we thus infer
  from (\ref{17.1}) and the nonegativity of $\varphi$ that
  \bas
	\int_0^\infty \io \frac{1}{(u+1)^2} |\nabla u|^2 \varphi
	&\le& \liminf_{\eps=\eps_j\searrow 0} \int_0^\infty \io \frac{1}{(\ue+1)^2} |\nabla \ue|^2 \varphi \\
	&=& - \int_0^\infty \io \ln (u+1)\varphi_t
	- \io \ln (u_0+1)\varphi(\cdot,0) \\
	& & - \int_0^\infty \io \ln (u+1)\Delta \varphi \\
	& & + \int_0^\infty \io \frac{u}{(u+1)^2} \nabla u \cdot \Big( S(x,u,v)\cdot \nabla v\Big) \cdot \varphi \\
	& & - \int_0^\infty \io \frac{u}{u+1} \Big(S(x,u,v)\cdot\nabla v\Big) \cdot \nabla \varphi
  \eas
  for any such test function $\varphi$, meaning that $u$ indeed is a global very weak $\phi$-supersolution of (\ref{14.1}).
\qed
Thereby our main result on global existence has actually been established already:\abs
\proofc of Theorem \ref{theo20}. \quad
  We only need to combine Lemma \ref{lem18} and Lemma \ref{lem17} with Lemma \ref{lem11}.
\qed
\mysection{Appendix}
Let us briefly collect some basic facts on approximation properties of the Steklov averages defined by
\bas
	(A_h w)(x,t) := \frac{1}{h} \int_{t-h}^t w(x,s)ds,
	\qquad x\in\Omega, \ t\in (0,T), \ h\in (0,1),
\eas
of a given function $w\in L^1(\Omega\times (-1,T)), T>0$.
\begin{lem}\label{lem_a1}
  As $h\searrow 0$, we have $A_h w \to w$ a.e.~in $\Omega \times (0,T)$.
\end{lem}
\proof
  Since $w\in L^1(\Omega \times (-1,T))$, there exists a null set $N\subset\Omega$ such that
  $(-1,T)\ni t \mapsto w(x,t)$ belongs to $L^1((-1,T))$ for all $x\in \Omega\setminus N$. Then for fixed 
  $x\in\Omega\setminus N$, by a known result one can pick a null set $N(x) \subset (0,T)$ such that each 
  $t\in (0,T)\setminus N(x)$ is a Lebesgue point of $(0,T) \ni \tilde t \mapsto w(x,\tilde t)$. This means that with
  \bas
	N_\star:=\Big\{ (x,t) \in \Omega \times (0,T) \ \Big| \ 
	x\in N \mbox{ or } (x\in \Omega\setminus N \mbox{ and } t\in N(x)) \Big\},
  \eas
  we have $(A_h w)(x,t) \to w(x,t)$ as $h\searrow 0$ for any $(x,t)\in (\Omega \times (0,T)) \setminus N_\star$.
  But since
  \bas
	\io \int_0^T \chi_{N_\star}(x,t)dtdx 
	&=& \int_N \int_0^T dtdx + \int_{\Omega\setminus N} \int_{N(x)} dtdx \\
	=|N| \cdot T + \int_{\Omega \setminus N} |N(x)|dx =0,
  \eas
  the Tonelli theorem ensures that $N_\star$ is a null set in $\Omega \times (0,T)$.
\qed
\begin{lem}\label{lem_a2}
  a) \ If $w\in L^p(\Omega \times (-1,T))$ for some $p\in (1,\infty)$, then $A_h w \wto w$ in $L^p(\Omega \times (0,T))$ as
  $h\searrow 0$. \\[1mm]
  b) \ If $w\in L^\infty(\Omega \times (-1,T))$, then $A_h w \wsto w$ in $L^\infty(\Omega\times (0,T))$ as $h\searrow 0$.
\end{lem}
\proof
  In view of a standard argument involving appropriate extraction of subsequences, Lemma \ref{lem_a1} and Egorov's theorem,
  it is sufficient to assert boundedness of $(A_h w)_{h\in (0,1)}$ in the respective spaces. 
  In the case in a) this follows on applying the H\"older inequality and Fubini's theorem in estimating
  \bas
	\|A_h w\|_{L^p(\Omega \times (0,T))}^p
	&=& \frac{1}{h^p} \io \int_0^T \bigg| \int_{t-h}^t w(x,s)ds \bigg|^p dtdx \\
	&\le& \frac{1}{h^p} \cdot h^{p-1} \io \int_0^T \int_{t-h}^t |w(x,s)|^p dsdtdx \\
	&=& \frac{1}{h} \io \int_0^T \int_s^{s+h} |w(x,s)|^p dtdsdx \\
	&=& \|w\|_{L^p(\Omega \times (-h,T-h))}^p 
	\le \|w\|_{L^p(\Omega \times (-1,T))}^p
	\qquad \mbox{for all } h\in (0,1).
  \eas
  In the situation in b) it is immediate that 
  $\|A_h w\|_{L^\infty(\Omega \times (0,T))} \le \|w\|_{L^\infty(\Omega \times (-1,T))}$ for all $h\in (0,1)$.
\qed
The derivation of the following criterion for strong convergence in $L^1$ is quite straightforward. Since we could not
find a precise reference in the literature, we include a short proof.
\begin{lem}\label{lem8}
  Let $N\ge 1$ and $M\subset \R^N$ be measurable, and suppose that $(w_j)_{j\in\N} \subset L^1(M)$ is such that $w_j\ge 0$
  a.e.~in $M$ for all $j\in\N$ and
  \be{8.1}
	w_j \wto w \quad \mbox{in } L^1(M)
	\qquad \mbox{and} \qquad
	w_j \to w \quad \mbox{a.e.~in } M
	\qquad \mbox{as } j\to\infty
  \ee
  with some $w\in L^1(M)$. Then
  \be{8.2}
	w_j \to w \quad \mbox{in } L^1(M)
	\qquad \mbox{as } j\to\infty
  \ee
\end{lem}
\proof
  Since $(w_j)_{j\in\N}$ is necessarily bounded in $L^1(M)$, the sequence $(\sqrt{w_j})_{j\in\N}$ is bounded in the 
  Hilbert space $L^2(M)$ and hence relatively compact in this space with respect to the weak topology. In view
  of the second assumption in (\ref{8.1}) and Egorov's theorem we thus have $\sqrt{w_j} \wto \sqrt{w}$ in $L^2(M)$
  as $j\to\infty$.\\
  Using constant test functions in the weak convergence statement in (\ref{8.1}), we moreover obtain that
  \bas
	\int_M (\sqrt{w_j})^2 = \int_M w_j\cdot 1 \to \int_M w\cdot 1 = \int_M (\sqrt{w})^2
	\qquad \mbox{as } j\to\infty,
  \eas
  which combined with the former yields $\sqrt{w_j} \to \sqrt{w}$ in $L^2(M)$ as $j\to\infty$.
  By means of the Cauchy-Schwarz inequality, this in turn implies that
  \bas
	\|w_j-w\|_{L^1(M)}
	&=& \int_M \Big|(\sqrt{w_j})^2-(\sqrt{w})^2\Big|
	= \int_M |\sqrt{w_j}-\sqrt{w}| \cdot |\sqrt{w_j}+\sqrt{w}| \\
	&\le& \Big\| \sqrt{w_j} - \sqrt{w} \Big\|_{L^2(M)} \cdot \Big( \|\sqrt{w_j}\|_{L^2(M)} + \|\sqrt{w}\|_{L^2(M)}\Big) \\
	&\to& 0
	\qquad \mbox{as } j\to\infty,
  \eas
  as claimed.
\qed
We finally note a useful consequence of the dominated convergence theorem.
\begin{lem}\label{lem19}
  Let $N\ge 1$ and $M\subset \R^N$ be measurable, and suppose that $(w_j)_{j\in\N} \subset L^\infty(M)$ and 
  $(z_j)_{j\in\N} \subset L^2(M)$ are such that
  \be{19.3}
	|w_j| \le C \quad \mbox{in } M \qquad \mbox{for all } j\in\N 
  \ee
  as well as
  \be{19.1}
	w_j \to w
	\qquad \mbox{a.e.~in } M
  \ee
  and
  \be{19.2}
	z_j \to z
	\qquad \mbox{in } L^2(M)
  \ee
  as $j\to\infty$ for some $C>0$, $w\in L^\infty(M)$ and $z\in L^2(M)$.
  Then 
  \be{19.4}
	w_j z_j \to wz
	\quad \mbox{in } L^2(M)
	\qquad \mbox{as } j\to\infty.
  \ee
\end{lem}
\proof
  We directly estimate
  \be{19.5}
	\int_N|w_j z_j - wz|^2 \le 2 \int_M (w_j-w)^2 z^2 + 2 \int_M w_j^2 (z_j-z)^2,
  \ee
  where by (\ref{19.3})-(\ref{19.2}) and the dominated convergence theorem we have
  \bas
	2 \int_M (w_j-w)^2 z^2 \to 0
	\qquad \mbox{as } j\to\infty.
  \eas
  Since thanks to (\ref{19.3}) and (\ref{19.2}) we know that also
  \bas
	2\int_M w_j^2(z_j-z)^2
	\le 2C^2 \|z_j-z\|_{L^2(M)}^2 \to 0
	\qquad \mbox{as } j\to\infty,
  \eas
  (\ref{19.5}) implies (\ref{19.4}).
\qed
{\bf Acknowledgement.} \quad
The author would like to thank Xinru Cao, Johannes Lankeit and the anonymous reviewer
for numerous fruitful comments on the topic of this paper.
\end{document}